\documentstyle[12pt]{article}

\font\twlgot =eufm10 scaled \magstep1
\font\egtgot =eufm8
\font\sevgot =eufm7
\font\twlmsb =msbm10 scaled \magstep1
\font\egtmsb =msbm8
\font\sevmsb =msbm7

\newfam\gotfam
\def\pgot{\fam\gotfam\twlgot}
\textfont\gotfam\twlgot
\scriptfont\gotfam\egtgot
\scriptscriptfont\gotfam\sevgot
\def\got{\protect\pgot}
\newfam\msbfam
\textfont\msbfam\twlmsb
\scriptfont\msbfam\egtmsb
\scriptscriptfont\msbfam\sevmsb
\def\Bbb{\protect\pBbb}
\def\pBbb{\relax\ifmmode\expandafter\Bb\else\typeout{You cann't use
Bbb in text mode}\fi}
\def\Bb #1{{\fam\msbfam\relax#1}}

\newcommand{\gF}{{\got F}}

\def\thebibliography#1{\bigskip\section*{}\bigskip\list
{$^{\arabic{enumi}}$}{\settowidth\labelwidth{#1}\leftmargin\labelwidth
\advance\leftmargin\labelsep
\usecounter{enumi}}
\def\newblock{\hskip .11em plus .33em minus .07em}
\sloppy\clubpenalty4000\widowpenalty4000
\sfcode`\.=1000\relax}

\let\Large=\large

\def\op#1{\mathop{\fam0 #1}\limits}

\newcommand{\id}{{\rm Id\,}}

\newcommand{\beq}{\begin{equation}}
\newcommand{\eeq}{\end{equation}}
\newcommand{\ben}{\begin{eqnarray}}
\newcommand{\een}{\end{eqnarray}}
\newcommand{\be}{\begin{eqnarray*}}
\newcommand{\ee}{\end{eqnarray*}}
\newcommand{\bea}{\begin{eqalph}}
\newcommand{\eea}{\end{eqalph}}
\newcommand{\cA}{{\cal A}}

\newcommand{\cV}{{\cal V}}

\newcommand{\cH}{{\cal H}}
\newcommand{\cF}{{\cal F}}
\newcommand{\cS}{{\cal S}}

\newcommand{\al}{\alpha}

\newcommand{\bt}{\beta}
\newcommand{\dl}{\delta}
\newcommand{\la}{\lambda}

\newcommand{\f}{\phi}
\newcommand{\om}{\omega}
\newcommand{\Om}{\Omega}
\newcommand{\m}{\mu}
\newcommand{\n}{\nu}
\newcommand{\g}{\gamma}
\newcommand{\G}{\Gamma}

\newcommand{\vt}{\vartheta}
\newcommand{\vf}{\varphi}

\newcommand{\di}{{\rm dim\,}}

\newcommand{\si}{\sigma}
\newcommand{\Si}{\Sigma}
\newcommand{\w}{\wedge}
\newcommand{\wt}{\widetilde}

\newcommand{\ol}{\overline}
\newcommand{\dr}{\partial}
\newcommand{\ar}{\op\longrightarrow}

\newcommand{\fl}{\flat}
\newcommand{\sh}{\sharp}

\let\ssection=\section
\renewcommand{\section}{\setcounter{equation}{0}\ssection}

\newcounter{eqalph}
\newcounter{equationa}
\newcounter{remark}
\newcounter{example}
\newcounter{theorem}
\newcounter{proposition}
\newcounter{lemma}
\newcounter{corollary}
\newcounter{definition}
\newenvironment{eqalph}{\stepcounter{equation}
\setcounter{equationa}{\value{equation}}
\setcounter{equation}{0}

\begin{eqnarray}}{\end{eqnarray}\setcounter{equation}{\value{equationa}}}

\def\therexample{\arabic{remark}}
\def\thetheorem{\arabic{theorem}}

\def\thedefinition{\arabic{definition}}

\newenvironment{proof}{\medskip\noindent
{\it Proof:}}{\medskip}

\newenvironment{theo}{\refstepcounter{theorem} \medskip
\noindent{\it Theorem \thetheorem:}}{\medskip}
\newenvironment{prop}{\refstepcounter{theorem} \medskip
\noindent{\it Proposition \thetheorem:}}{\medskip}
\newenvironment{lem}{\refstepcounter{theorem} \medskip
\noindent{\it Lemma \thetheorem:}}{\medskip}

\newenvironment{defi}{\refstepcounter{definition} \medskip
\noindent{\it Definition \thedefinition:}}{\medskip}

\newcommand{\mar}[1]{}

\begin{document}
\hbox{}

{\parindent=0pt

{\Large \bf Bi-Hamiltonian partially integrable
systems}

\bigskip

{\sc G.Giachetta}\footnote{Electronic mail: giovanni.giachetta@unicam.it},

{\sl Department of Mathematics and Informatics, University of Camerino,
62032 Camerino (MC), Italy}

\medskip

{\sc L.Mangiarotti}\footnote{Electronic mail: luigi.mangiarotti@unicam.it},

{\sl Department of Mathematics and Informatics, University of Camerino,
62032 Camerino (MC), Italy}

\medskip

{\sc G. Sardanashvily}\footnote{Electronic mail:
sard@grav.phys.msu.su}

{\sl Department of Theoretical Physics,
Moscow State University, 117234 Moscow, Russia}

\bigskip
Given a first order dynamical system possessing a
commutative algebra of dynamical
symmetries, we show that,
under certain conditions, there exists
a Poisson structure
on an open neighbourhood of its regular
(not necessarily compact) invariant manifold
which makes this dynamical system into a
partially integrable Hamiltonian system.
This Poisson structure is by no means unique.
Bi-Hamiltonian partially integrable systems
are described in some detail. As an outcome, we state the
conditions of
quasi-periodic stability (the KAM theorem)
for partially integrable Hamiltonian systems.
}

\bigskip
\bigskip

\noindent
{\bf I. INTRODUCTION}
\bigskip

Given a smooth real manifold $Z$, let we have $k$
mutually commutative vector fields $\{\vt_\la\}$ which are independent
almost everywhere on $Z$, i.e., the set of points where
the multivector field $\op\w^k\vt_\la$ vanishes is nowhere dense.
We denote by $\cS\subset C^\infty(Z)$ the $\Bbb R$-subring of smooth real
functions $f$ on $Z$
whose derivations $\vt_\la\rfloor df$ vanish for all $\vt_\la$.
Let $\cA$ be
the $k$-dimensional $\cS$-Lie algebra generated by
the vector fields $\{\vt_\la\}$.
One can think of
one of its elements as being a first order dynamic equation
on $Z$ and of the other as being the dynamical
symmetries. Accordingly, elements of $\cS$ are regarded as
integrals of motion. For the sake of brevity, we agree to call
$\cA$ a dynamical algebra.

Completely and partially integrable systems
on symplectic manifolds$^1$
and broadly integrable
dynamical systems of Bogoyavlenkij$^{2,3}$
exemplify finite-dimensional commutative dynamical
algebras.
Recall that, given a symplectic manifold $(Z,\Om)$, we have a 
partially integrable
system (henceforth PIS) if there exist $1\leq k\leq \di Z/2$
smooth real functions $\{H_\la\}$
in involution which are independent almost
everywhere on $Z$, i.e., the set of points where the $k$-form $\op\w^k dH_\la$
vanishes is nowhere dense.
The Hamiltonian vector fields $\vt_\la$
of functions $H_\la$ mutually commute and are independent almost
everywhere. They make up a commutative dynamical algebra
over the  Poisson subalgebra $\cS$ of elements of $C^\infty(Z)$ commuting
with all the functions $H_\la$.

An important peculiarity of a
finite-dimensional commutative dynamical algebra $\cA$ is that
its regular invariant manifolds are toroidal cylinders
$\Bbb R^{k-m}\times T^m$.
At the same time, no preferable Poisson structure is associated to
a commutative Lie algebra $\cA$ because its Lie--Poisson structure
is zero. Therefore, we are free with analyzing
different Poisson structures which make $\cA$ into a Hamiltonian system.
However, this analysis essentially differs from that of
noncommutative integrable systems (see Ref. [4] for a survey).
One has investigated different symplectic structures around
invariant tori of commutative integrable systems.$^{2,3,5,6}$
For instance, the classical Liouville--Arnold theorem$^{1,7}$
and the Nekhoroshev theorem$^{8,9}$ state that, under certain
conditions, every symplectic structure making a commutative dynamical 
algebra into
a Hamiltonian system takes a canonical form around a compact 
invariant manifold.

Our goal is to describe all Poisson structures bringing a commutative 
dynamical algebra into
a PIS near its regular invariant manifold, which
need not be compact.

\begin{defi} \label{cc1} \mar{cc1}
A $k$-dimensional commutative dynamical algebra on a regular Poisson 
manifold $(Z,w)$
is said to be a PIS if:
(a) $\cA$ is
generated by Hamiltonian vector fields of $k$ almost everywhere 
independent integrals
of motion $H_\la\in C^\infty(Z)$ in involution;
(b) all elements of $\cS\subset C^\infty(Z)$ are mutually in involution.
\end{defi}

It follows at once from this definition that the Poisson structure 
$w$ is at least of rank $2k$
and $\cS$ is a commutative Poisson algebra. If $2k=\di Z$, we have
a completely integrable system on a symplectic manifold.

Theorem \ref{bi0} below states that:

(i) under certain conditions,
an open neighbourhood $U$
of a regular invariant manifold $M$ of the dynamical system $\cA$ is
a trivial
principal bundle
\mar{z10}\beq
U=N\times(\Bbb R^{k-m}\times T^m)\ar^\pi N \label{z10}
\eeq
over a domain $N\subset \Bbb R^{\di Z-k}$ with the
structure group $\Bbb R^{k-m}\times T^m$;

(ii) the toroidal domain (\ref{z10}) is provided
with a Poisson structure $w$
such that $(w,\cA)$ is a PIS in accordance with Definition \ref{cc1}.

Note that a trivial fibration in invariant manifolds is a standard property
of integrable systems.$^{1-3,8-12}$
However, there exists a well-known
obstruction to its global extension in the case of compact invariant 
manifolds,$^{13}$
and there is an additional obstruction
similar to that in Ref. [3] if invariant manifolds are noncompact.

A Poisson structure in  Theorem \ref{bi0}
is by no means unique.
Let the toroidal domain (\ref{z10})
be provided with bundle coordinates
$(r^A,y^\la)$, where $(r^A)$ are coordinates on $N$ and $(y^\la)=(t^a,\vf^i)$
are standard coordinates on the toroidal cylinder $\Bbb R^{k-m}\times T^m$.
It is readily observed that, if a Poisson bivector field
on the toroidal domain $U$ satisfies Definition \ref{cc1}, it takes the
form
\mar{bi20}\beq
w=w_1+w_2=w^{A\la}(r^B)\dr_A\w\dr_\la + w^{\m\nu}(r^B,y^\la)\dr_\m\w \dr_\nu.
\label{bi20}
\eeq
The converse also holds (see Theorem \ref{bi92} below).
For any Poisson bivector field $w$ (\ref{bi20}) of rank
$2k$ on $U$,
there exists a toroidal domain $U'\subset U$
such that $(w,\cA)$ is a PIS on $U'$.
Moreover, Theorem \ref{bi100} in Section III states that
there is a toroidal domain
  $U'$
such that, restricted to $U'$, this Poisson bivector field
takes the canonical form (\ref{bi101}).

Now, let $w$ and $w'$ be two different Poisson structures (\ref{bi20})
which make a commutative dynamical algebra $\cA$ into the different
PISs $(w,\cA)$ and $(w',\cA)$.

\begin{defi} \label{cc2} \mar{cc2}
We agree to call the triple $(w,w',\cA)$
a bi-Hamiltonian PIS if
any Hamiltonian vector field $\vt\in\cA$ with respect to
$w$ possesses the same
  Hamilton representation
\mar{bi71}\beq
\vt=-w\lfloor df=-w'\lfloor df, \qquad f\in\cS, \label{bi71}
\eeq
relative to $w'$, and {\it vice versa}.
\end{defi}

Definition \ref{cc2} establishes a {\it sui generis}
equivalence between the PISs
$(w,\cA)$ and $(w',\cA)$. Theorem \ref{bi72}
below states that the triple $(w,w',\cA)$ is
a bi-Hamiltonian PIS in accordance
with Definition \ref{cc2} iff
the  Poisson bivector fields $w$ and $w'$
(\ref{bi20}) differ only
in the second terms $w_2$ and $w'_2$. Moreover, these
Poisson bivector fields admit a recursion operator.

Let now $\cA$ be a commutative dynamical algebra
associated to a PIS
on a symplectic manifold $(Z,\Om)$. In this case,
condition (b) in Definition \ref{cc1} is not necessarily satisfied, unless
it is a completely integrable system.
Nevertheless, there exists
a Poisson structure $w$ of rank $2k$
on the toroidal domain (\ref{z10}) such that,
with respect to $w$,
all integrals of motion $H_\la$ of the original PIS remain
to be in involution, and they possess the same Hamiltonian vector
fields $\vt_\la$ (see Theorem \ref{dd3} below).
Therefore, one can think of
the triple $(\Om,w,\{H_\la\})$ as being
a special bi-Hamiltonian system, though it fails to satisfy 
Definition \ref{cc2}.
Conversely, if $Z$ is even-dimensional,
any Poisson bivector field $w$ (\ref{bi20}) setting a PIS $(w,\cA)$
is extended
to an appropriate symplectic structure $\Om$ such that $(\Om,\cA)$ is a PIS
on the symplectic manifold $(Z,\Om)$.

There are several reasons in order to make a commutative dynamical algebra
$\cA$ into a Hamiltonian
one. For instance, one can quantize $\cA$ around its invariant manifold
by quantizing the Poisson algebra $\cS$.$^{12,14}$ Of course, 
quantization of $\cA$
with respect to different Poisson structures need not be equivalent.
However, we focus on another result. In Section V, we show that, introducing
an appropriate Poisson structure and using the methods in Ref. [10], one
can extend the well-known KAM theorem to PISs.

\bigskip
\bigskip

\noindent
{\bf II. SEMILOCAL GEOMETRY AROUND AN INVARIANT MANIFOLD}
\bigskip

Given a $k$-dimensional commutative dynamical algebra $\cA$
on a smooth manifold $Z$,
let $\cV$ be the smooth involutive distribution on $Z$
spanned by the vector fields $\{\vt_\la\}$, and let
$G$ be the group of local diffeomorphisms
of $Z$ generated by the flows of these vector fields
(we follow the terminology of Ref. [15]).
Maximal integral manifolds of $\cV$
are the orbits of $G$, and are
invariant manifolds of $\cA$.$^{15}$
Let $z\in Z$ be a regular
point of the distribution $\cV$, i.e., $\op\w^k\vt_\la(z)\neq 0$.
Since the group $G$ preserves $\op\w^k\vt_\la$, the maximal integral manifold
$M$ of $\cV$ through $z$ is also regular. Furthermore, there exists an open
neighbourhood $U$ of $M$ such that, restricted to $U$, the distribution
$\cV$ is regular and yields a foliation $\gF$ of $U$.

\begin{theo} \label{bi0} \mar{bi0}
Let us suppose that: (i)
the vector fields $\vt_\la$ on $U$ are complete, (ii)
the foliation $\gF$ of $U$ admits a transversal manifold $\Si$
and its holonomy pseudogroup on $\Si$ is trivial, (iii) the leaves of this
foliation are mutually diffeomorphic.
Then the following hold.

(I) There exists an open neighbourhood of $M$, say $U$ again, which 
is the trivial
principal bundle (\ref{z10})
over a domain $N\subset \Bbb R^{\di Z-k}$ with the
structure group $\Bbb R^{k-m}\times T^m$.

(II) If $2k\leq\di Z$, there exists a Poisson structure $w$
of rank $2k$ on $U$ such that $(w,\cA)$ is a PIS in accordance
with Definition \ref{cc1}.
\end{theo}

Let us note the following. Condition (i)
states that $G$ is a group of diffeomorphisms of $U$.
Condition (ii) is equivalent to the assumption that
$U\to U/G$ is a fibered manifold.$^{16}$
Each fiber $M_r$, $r\in N$, of this fibered manifold
admits a free transitive
action of the group $G_r=G/K_r$, where $K_r$ is the
isotropy group of an arbitrary point of $M_r$. In accordance with 
condition (iii),
all the groups $G_r$, $r\in N$, are isomorphic
to the toroidal  cylinder group $\Bbb R^{k-m}\times T^m$ for some 
$0\leq m\leq k$.
The goal is to define these isomorphisms so that
they provide a smooth action of $\Bbb R^{k-m}\times T^m$ in $U$.
We follow the proof in Refs. [7,17] generalized to noncompact
invariant manifolds.
We establish a particular
trivialization (\ref{z10}) such that the generators $\vt_\la$ of the algebra
$\cA$ take the specific form (\ref{ww25}).
Part (II) of Theorem \ref{bi0} is based on this trivialization

\begin{proof} (I).
By virtue of the condition (ii), the foliation $\gF$ of $U$ is a 
fibered manifold
\mar{d20}\beq
\pi:U\to N, \label{d20}
\eeq
admitting a section $\si$ such that and $\Si=\si(N)$.$^{16}$
Since the vector fields $\vt_\la$ on $U$
are complete and mutually commutative, the group $G$
of their flows is an additive Lie group
of diffeomorphism of $U$. Its group space is a
vector space $\Bbb R^k$ coordinated by parameters $(s^\la)$
of the flows with respect to the basis $\{e_\la=\vt_\la\}$.
Since vector fields $\vt_\la$ are independent everywhere on $U$, the 
action of $\Bbb R^k$
in $U$ is locally free, i.e., isotropy groups of points of $U$ are
discrete subgroups of the group $\Bbb R^k$. Its orbits
are fibers  of the fibered manifold (\ref{d20}).
Given a point $r\in N$, the action of $\Bbb R^k$
in the fiber $M_r=\pi^{-1}(r)$ factorizes as
\mar{d4}\beq
\Bbb R^k\times M_r\to G_r\times M_r\to M_r \label{d4}
\eeq
through the free transitive
action in $M_r$ of the factor group $G_r=\Bbb R^k/K_r$, where $K_r$ 
is the isotropy group of
an arbitrary point of $M_r$. It is the same group for all points of $M_r$
because $\Bbb R^k$ is a commutative group.
Since the fibers $M_r$ are mutually diffeomorphic, all isotropy groups $K_r$
are isomorphic to the group $\Bbb Z^m$ for some fixed $0\leq m\leq k$.
Accordingly, the groups
$G_r$ are isomorphic to the additive group $\Bbb R^{k-m}\times T^m$.
Let us bring the fibered manifold $U\to N$ (\ref{d20}) into
a principal bundle with
the structure group $G_0$, where we denote $\{0\}=\pi(M)$. For this purpose,
let us determine
isomorphisms $\rho_r: G_0\to G_r$ of the group $G_0$ to the groups 
$G_r$, $r\in N$.
Then, a desired fiberwise action of $G_0$ in $U$ is defined by the law
\mar{d5}\beq
G_0\times M_r\to\rho_r(G_0)\times M_r\to M_r. \label{d5}
\eeq
Generators of each isotropy subgroup
$K_r$ of $\Bbb R^k$ are given by $m$ linearly independent vectors of the group
space $\Bbb R^k$. One can show that there exist ordered collections 
of generators
  $(v_1(r),\ldots,v_m(r))$
of the groups $K_r$ such that $r\mapsto v_i(r)$
are smooth $\Bbb R^k$-valued fields on $N$. Indeed, given a vector $v_i(0)$
and a section $\si$ of the fibered manifold (\ref{d20}),
each field $v_i(r)=(s^\al(r))$ is the unique smooth solution of the equation
\be
g(s^\al)\si(r)=\si(r), \qquad  (s^\al(0))=v_i(0),
\ee
on an open neighbourhood of $\{0\}$.
Let us consider the decomposition
\be
v_i(0)=B_i^a(0) e_a + C_i^j(0) e_j, \qquad a=1,\ldots,k-m, \qquad 
j=1,\ldots, m,
\ee
where $C_i^j(0)$ is a nondegenerate matrix.
Since the fields $v_i(r)$ are smooth, there exists an open 
neighbourhood of $\{0\}$,
say $N$ again, where the matrices $C_i^j(r)$ are nondegenerate. Then,
\mar{d6}\beq
A_r=\left(
\begin{array}{ccc}
\id & \qquad & (B(r)-B(0))C^{-1}(0) \\
0 & & C(r)C^{-1}(0)
\end{array}
\right) \label{d6}
\eeq
is
a unique linear
morphism of the vector space $\Bbb R^k$ which transforms
the frame $v_\la(0)=\{e_a,v_i(0)\}$
into the frame $v_\la(r)=\{e_a,v_i(r)\}$.
Since it is also an automorphism of the group $\Bbb R^k$
sending $K_0$ onto $K_r$, we obtain a desired isomorphism
$\rho_r$ of the group $G_0$ to the group $G_r$. Let an element $g$ of 
the group $G_0$
be the coset of an element $g(s^\la)$ of the group $\Bbb R^k$. Then, it
acts in $M_r$ by the rule (\ref{d5}) just as the element 
$g((A_r^{-1})^\la_\bt s^\bt)$
of the group $\Bbb R^k$ does. Since entries of the matrix $A$ (\ref{d6}) are
smooth functions on $N$, this action of the group $G_0$ in $U$ is 
smooth. It is free, and
$U/G_0=N$. Then, the fibered manifold $U\to N$ is a trivial principal 
bundle with the structure
group $G_0$.
Given a section $\si$ of the principal bundle $U\to N$, its 
trivialization $U=N\times G_0$
is defined by assigning the points $\rho^{-1}(g_r)$ of the
group space $G_0$ to the points
$g_r\si(r)$, $g_r\in G_r$, of a fiber $M_r$.
Let us endow $G_0$ with the standard
coordinate atlas
$(y^\la)=(t^a,\vf^i)$ of the group $\Bbb R^{k-m}\times T^m$. Then, we provide
$U$ with the trivialization (\ref{z10})
with respect to the coordinates $(r^A,t^a,\vf^i)$,
where $(r^A)$, $A=1,\ldots, \di Z-k,$ are coordinates on the base $N$.
The vector fields $\vt_\la$ on $U$ relative to these coordinates read
\mar{ww25}\beq
\vt_a=\dr_a, \qquad \vt_i=-(BC^{-1})^a_i(r)\dr_a +
(C^{-1})_i^k(r)\dr_k.\label{ww25}
\eeq
Accordingly,
the subring $\cS$ restricted to $U$ is
the pull-back $\pi^*C^\infty(N)$ onto $U$ of the ring
of smooth functions on $N$.

(II). Let us split the coordinates $(r^A)$ into some
$k$ coordinates $(I_\la)$ and $\di Z- 2k$ coordinates
$(z^A)$. Then, we can provide the toroidal domain
$U$ (\ref{z10}) with the Poisson bivector field
\mar{d26'}\beq
w=\dr^\la\w\dr_\la \label{d26'}
\eeq
of rank $2k$. The independent complete vector fields
$\dr_a$ and $\dr_i$ are Hamiltonian vector fields of the functions
$H_a=I_a$ and $H_i=I_i$ on $U$ which are in involution with respect to
the Poisson bracket
\mar{bi12}\beq
\{f,f'\}=\dr^\la f\dr_\la f'-\dr_\la f\dr^\la f' \label{bi12}
\eeq
defined by the bivector field (\ref{d26'}).
By virtue of the
expression (\ref{ww25}), the Hamiltonian vector fields $\{\dr_\la\}$
generate the $\cS$-algebra $\cA$.
\end{proof}

\bigskip
\bigskip

\noindent
{\bf III. POISSON STRUCTURES AROUND AN INVARIANT MANIFOLD}
\bigskip

\begin{theo} \label{bi92} \mar{bi92}
For any Poisson bivector field $w$ (\ref{bi20}) of rank
$2k$ on $U$,
there exists a toroidal domain $U'\subset U$
such that $(w,\cA)$ is a PIS on $U'$.
\end{theo}

It is readily observed that any Poisson bivector
field $w$ (\ref{bi20}) fulfills condition (b) in Definition \ref{cc1},
but condition (a) imposes a restriction on the
toroidal domain $U$.
The key point is that the characteristic foliation $\cF$
of $U$ yielded by the Poisson bivector fields $w$ (\ref{bi20})
is the pull-back of a $k$-dimensional
foliation $\cF_N$ of the base $N$, which is defined by the
first summand $w_1$ (\ref{bi20}) of $w$.
With respect to the adapted coordinates $(J_\la,z^A)$, $\la=1,\ldots, k$,
on the foliated
manifold $(N,\cF_N)$, the
Poisson bivector field $w$ reads
\mar{bi42}\beq
w= w^\m_\n(J_\la,z^A)\dr^\n\w \dr_\m +
w^{\m\n}(J_\la,z^A,y^\la)\dr_\m\w \dr_\n. \label{bi42}
\eeq
Then, condition (a) in Definition \ref{cc1} is
satisfied if $N'\subset N$ is a domain of a coordinate chart $(J_\la,z^A)$
of the foliation $\cF_N$. In this case, the
dynamical algebra $\cA$ on the toroidal domain $U'=\pi^{-1}(N')$ is 
generated by the
Hamiltonian vector fields
\mar{bi93}\beq
\vt_\la=-w\lfloor dJ_\la=w^\m_\la\dr_\m \label{bi93}
\eeq
of the $k$ independent functions $H_\la=J_\la$.

\begin{proof}
The characteristic distribution of the Poisson bivector field $w$ 
(\ref{bi20}) is spanned by
the Hamiltonian vector fields
\mar{bi21}\beq
v^A=-w\lfloor dr^A=w^{A\m}\dr_\m \label{bi21}
\eeq
and the vector fields
\be
w\lfloor dy^\la= w^{A\la}\dr_A + 2w^{\m\la}\dr_\m.
\ee
Since $w$ is of rank $2k$, the vector fields $\dr_\m$ can be 
expressed into the vector
fields $v^A$ (\ref{bi21}).
Hence, the characteristic distribution of $w$ is spanned by
the Hamiltonian vector fields $v^A$ (\ref{bi21}) and the vector fields
\mar{bi25}\beq
v^\la=w^{A\la}\dr_A. \label{bi25}
\eeq
The vector fields (\ref{bi25}) are projected onto $N$. Moreover, one 
can derive from the relation
$[w,w]=0$ that they
generate a Lie algebra and, consequently, span an involutive
distribution $\cV_N$ of rank $k$ on $N$.
Let $\cF_N$ denote the corresponding foliation of $N$. We consider the
pull-back $\cF=\pi^*\cF_N$ of this foliation onto $U$ by the trivial
fibration $\pi$.$^{16}$ Its leaves are the inverse images 
$\pi^{-1}(F_N)$ of leaves $F_N$
of the foliation $\cF_N$, and so is its characteristic
distribution $T\cF=(T\pi)^{-1}(\cV_N)$.
This distribution is spanned by the vector fields $v^\la$ (\ref{bi25})
on $U$ and the vertical vector fields
on $U\to N$, namely,  the vector fields $v^A$ (\ref{bi21}) generating
the algebra $\cA$. Hence, $T\cF$ is the characteristic
distribution of the Poisson bivector field $w$. Furthermore, since $U\to N$
  is a trivial bundle,
each leaf $\pi^{-1}(F_N)$ of the pull-back foliation $\cF$ is the 
manifold product
of a leaf $F_N$ of $N$ and the toroidal cylinder $\Bbb R^{k-m}\times 
T^m$. It follows that
the foliated manifold $(U,\cF)$ can be provided with an adapted 
coordinate atlas
\mar{bi26}\beq
\{(U_\iota,J_\la,z^A,y^\la)\}, \qquad \la=1,\ldots, k, \qquad 
A=1,\ldots,\di Z-2k,
\label{bi26}
\eeq
  such that $(J_\la,z^A)$ are adapted coordinates on
the foliated manifold $(N,\cF_N)$, i.e., transition functions of coordinates
$z^A$ are independent of $J_\la$, while transition functions of coordinates
$(y^\la)=(t^a,\vf^\la)$ on the toroidal cylinder $\Bbb R^{k-m}\times T^m$
are independent of coordinates $J_\la$ and $z^A$.
With respect to these coordinates, the Poisson bivector
field (\ref{bi20}) takes the form (\ref{bi42}).
Let $N'$ be the domain of a coordinate chart (\ref{bi26}). Then, the
dynamical algebra $\cA$ on the toroidal domain $U'=\pi^{-1}(N')$ is 
generated by the
Hamiltonian vector fields $\vt_\la$ (\ref{bi93}) of functions
$H_\la=J_\la$.
\end{proof}

Note that the coefficients $w^{\m\nu}$ in the expressions
(\ref{bi20}) and (\ref{bi42}) are affine in coordinates $y^\la$
because of the relation $[w,w]=0$ and, consequently, are constant
on tori. Furthermore, one can improve the expression (\ref{bi42}) as follows.

\begin{theo} \label{bi100} \mar{bi100}
Given a PIS
$(w,\cA)$ on a Poisson manifold $(w,U)$,
there exists a toroidal domain $U'\subset U$
equipped with partial action-angle coordinates
$(I_a,I_i,z^A, x^a,\f^i)$
such that, restricted to $U'$, a Poisson bivector field
takes the canonical
form
\mar{bi101}\beq
w=\dr^a\w \dr_a + \dr^i\w \dr_i,
\label{bi101}
\eeq
while the dynamical algebra $\cA$ is generated by Hamiltonian vector fields
of the action coordinate functions $H_a=I_a$, $H_i=I_i$.
\end{theo}

Theorem \ref{bi100}
extends the Liouville--Arnold theorem to the case of a Poisson structure
and a noncompact invariant manifold.
To prove it (see Appendix A), we reformulate the
proof of the Liouville--Arnold theorem for noncompact invariant manifolds
in Refs. [11,12] in terms of a leafwise
symplectic structure.

Given a dynamic equation $\xi\in\cA$,
it may happen
that no Poisson
bivector field (\ref{bi20}) makes $\xi$ into a Hamilton equation.
If $\xi$ is a nowhere
vanishing complete vector field whose trajectories are
not located in tori, one can choose $\xi$ as one of the generators,
e.g., $\xi=\vt_1$
in Theorem \ref{bi0} so that $U$ can be
provided with a trivialization such that
$\xi=\vt_1=\dr_1$ in the expression (\ref{ww25}). Then,
the Poisson structure (\ref{d26'}) brings $\xi$
into a Hamilton equation.
This improves the well-known result of
Hojman$^{18}$ that,
under certain conditions, a first order
dynamic equation can be brought into a Hamilton
one with respect to a Poisson structure of rank 2.
Moreover, any dynamic equation $\xi$
on $U$ gives rise to an
equivalent Hamilton equation
$\dr_t +\xi$ of time-dependent mechanics on
$\Bbb R^2\times U$.$^{11,12}$

\bigskip
\bigskip

\noindent
{\bf IV. BI-HAMILTONIAN STRUCTURES}
\bigskip

Now, let $w$ and $w'$ be two different Poisson structures (\ref{bi20})
on the toroidal domain (\ref{z10})
which make a commutative dynamical algebra $\cA$ into two different
PISs $(w,\cA)$ and $(w',\cA)$.

\begin{theo} \label{bi72} \mar{bi72}
(I) The triple $(w,w',\cA)$ is
a bi-Hamiltonian system PIS in accordance
with Definition \ref{cc2} iff
the  Poisson bivector fields $w$ and $w'$
(\ref{bi20}) differ only
in the second terms $w_2$ and $w'_2$. (II) These Poisson bivector 
fields admit a
recursion operator.
\end{theo}

\begin{proof} (I). It is easily justified that, if Poisson bivector 
fields $w$ (\ref{bi20}) fulfil
Definition \ref{cc2}, they are
distinguished only by the second summand $w_2$. Conversely,
as follows from the proof of Theorem \ref{bi92},
the characteristic distribution of a Poisson bivector field $w$ (\ref{bi20})
is spanned by the vector fields (\ref{bi21}) and (\ref{bi25}).
Hence, all Poisson bivector fields
$w$ (\ref{bi20}) distinguished only by the second summand $w_2$ have
the same characteristic distribution, and
they bring $\cA$ into a PIS on the same toroidal domain $U'$. Then,
the condition in Definition \ref{cc2} is easily justified.

(II). The result follows from forthcoming Lemma
\ref{p1}.
\end{proof}

Given  a smooth real manifold $X$, let
$w$ and $w'$ be Poisson bivector fields
of rank $2k$ on $X$, and let
\mar{p4}\beq
w^\sh: T^*X\to TX, \qquad w'^\sh: T^*X\to TX \label{p4}
\eeq
be the corresponding bundle homomorphisms. A
tangent-valued one-form $R$ on $X$
yields bundle endomorphisms
\mar{bi90}\beq
R: TX\to TX, \qquad R^*: T^*X\to T^*X. \label{bi90}
\eeq
It is called a recursion operator if
\mar{p0}\beq
w'^\sh=R\circ w^\sh=w^\sh\circ R^*. \label{p0}
\eeq
Given a Poisson bivector field $w$ and
a tangent valued one-form $R$ such that $R\circ w^\sh=w^\sh\circ R^*$,
the well-known sufficient condition for $R\circ w^\sh$ to be a Poisson bivector
field is that the Nijenhuis torsion of $R$ and the Magri--Morosi
concomitant of $R$ and $w$ vanish.$^{19,20}$
However, as we will see later, recursion operators between Poisson bivector
fields in Theorem \ref{bi72} need not satisfy these conditions.

\begin{lem} \label{p1} \mar{p1}
A recursion operator between Poisson structures
of the same rank exists iff their characteristic distributions
coincide.
\end{lem}

\begin{proof} It follows from the equalities
(\ref{p0}) that a recursion operator $R$
sends the characteristic distribution of $w$ to
that of $w'$, and these distributions coincide if
$w$ and $w'$ are of the same rank.
Conversely, let regular Poisson structures $w$ and $w'$
possess the same characteristic distribution
$T\cF\to TX$ tangent to a foliation $\cF$ of $X$.
Let $T\cF^*\to X$ be the dual of $T\cF\to X$, and let
\mar{p2,3}\ben
&& 0\to T\cF \ar^{i_\cF} TX \ar TX/T\cF\to 0, \label{p2} \\
&& 0\to {\rm Ann}\,T\cF\ar T^*X\ar^{i^*_\cF} T\cF^* \to 0, \label{p3}
\een
be the corresponding exact sequences.
The bundle homomorphisms $w^\sh$ and $w'^\sh$ (\ref{p4})
factorize in a unique fashion
\be
&& w^\sh:
T^*X\ar^{i^*_\cF} T\cF^*\ar^{w^\sh_\cF}
T\cF\ar^{i_\cF} TX, \\
&&  w'^\sh:
T^*X\ar^{i^*_\cF} T\cF^*\ar^{w'^\sh_\cF}
T\cF\ar^{i_\cF} TX
\ee
through the bundle isomorphisms
\be
w_\cF^\sh: T\cF^*\to T\cF,  \qquad w'^\sh_\cF: T\cF^*\to T\cF.
\ee
Let us consider the inverse isomorphisms
\mar{p13}\beq
w_\cF^\fl : T\cF\to T\cF^*, \qquad w'^\fl_\cF : T\cF\to T\cF^* \label{p13}
\eeq
and the compositions
\mar{p10}\beq
R_\cF= w'^\sh_\cF\circ w_\cF^\fl: T\cF\to T\cF, \qquad
R_\cF^*= w_\cF^\fl \circ w'^\sh_\cF: T\cF^*\to T\cF^*. \label{p10}
\eeq
There is the obvious relation
\be
w'^\sh_\cF=R_\cF\circ w^\sh_\cF=  w^\sh_\cF\circ R^*_\cF.
\ee
In order to obtain a recursion operator (\ref{p0}), it suffices
to extend the morphisms $R_\cF$ and $R_\cF^*$ (\ref{p10}) onto $TX$ 
and $T^*X$, respectively.
  For this purpose, let us consider a splitting
\be
\zeta: TX\to T\cF, \qquad TX=T\cF\oplus (\id-i_\cF\circ\zeta)TX=T\cF\oplus E,
\ee
of the exact sequence (\ref{p2}) and the dual splitting
\be
\zeta^* :T\cF^*\to T^*X, \qquad T^*X=\zeta^*(T\cF^*)\oplus
(\id-\zeta^*\circ i^*_\cF)T^*X= \zeta^*(T\cF^*)\oplus E'
\ee
of the exact sequence (\ref{p3}). Then, the desired extensions are
\be
R:=R_\cF\times \id E, \qquad R^*:=(\zeta^*\circ R^*_\cF)\times \id E'.
\ee
This recursion operator is invertible, i.e.,
the morphisms (\ref{bi90}) are bundle isomorphisms.
\end{proof}

For instance,  the Poisson bivector field $w$ (\ref{bi20}) and the
Poisson bivector field
\mar{bi22}\beq
w_0=w^{A\la}(r)\dr_A\w\dr_\la \label{bi22}
\eeq
admit a recursion operator $w^\sh_0=R\circ w^\sh$ whose entries are 
given by the equalities
\mar{bi24}\beq
R^A_B=\dl^A_B, \qquad R^\m_\n=\dl^\m_\n, \qquad R^A_\la=0, \qquad 
w^{\m\la}=R^\la_Bw^{B\m}.
\label{bi24}
\eeq
Its Nijenhuis torsion fails to vanish, unless coefficients $w^{\m\la}$
are independent of coordinates $y^\la$.

Turn now to the case of a commutative dynamical algebra $\cA$
defined by a PIS
on a symplectic manifold $(Z,\Om)$. The following generalization
of the Nekhoroshev
theorem to noncompact invariant manifolds addresses such a system.

\begin{theo} \mar{dd3} \label{dd3}
  Let $(\Om,\{H_\la\},\vt_\la)$
be a $k$-dimensional PIS on a $2n$-dimensional symplectic
manifold $(Z,\Om)$.
Let the distribution $\cV$, its regular integral manifold $M$,
an open neighbourhood $U$ of $M$, and the foliation $\gF$ of $U$ be as those
in Theorem \ref{bi0}.
Under conditions (i) -- (iii) of Theorem \ref{bi0}, the following hold.

(I) There exists an open neighbourhood of $M$, say $U$ again, which 
is the trivial
bundle (\ref{z10})
in toroidal cylinders $\Bbb R^{k-m}\times T^m$  over a domain 
$N\subset \Bbb R^{2n-k}$.

(II) It is provided with the partial action-angle coordinates
$(I_\la, z^A, y^\la)$
such that the functions $H_\la$ depend only on the
action coordinates $I_\la$ and the symplectic form $\Om$ on $U$
reads
\mar{d26}\beq
\Om= dI_\la\w d y^\la +\Om_{AB}(I_\m,z^C) dz^A\w dz^B +
\Om_A^\la(I_\m,z^C) dI_\la\w dz^A.
\label{d26}
\eeq

(III) There exists a Darboux coordinate
chart $Q\times \Bbb R^{k-m}\times T^m\subset U$, foliated in toroidal cylinders
$\Bbb R^{k-m}\times T^m$  and provided with
coordinates $(I_\la,p_s,q^s,\wt y^\la)$ such that the symplectic form 
$\Om$ (\ref{d26})
on this chart takes the canonical form
\mar{dd12}\beq
\Om= dI_\la\w d \wt y^\la + dp_s\w dq^s. \label{dd12}
\eeq
\end{theo}

This theorem is proved in Appendix B.
Part (I)
repeats exactly that of Theorem \ref{bi0}, while the proof of part (II)
follows that of Theorem \ref{bi100}.
The proof of part (III) is a generalization of that
of Proposition 1 in Ref. [21] to noncompact invariant manifolds.
As follows from the expression (\ref{dd12}), the PIS
in Theorem
\ref{dd3} can be extended to a completely integrable system on some open
neighbourhood of $M$, but Hamiltonian vector field of its additional
local integrals of motion fail to be complete.

A glance at the symplectic form $\Om$ (\ref{d26}) shows that there exists
a Poisson structure $w$ of rank $2k$, e.g., $w=\dr^\la\w \dr_\la$
on $U$ such that,
with respect to $w$,
the integrals of motion $H_\la$ of the original PIS remain
to be in involution, and they possess the same Hamiltonian vector
fields $\vt_\la$. Hence, $(\Om,w,\{H_\la\})$ is the above mentioned
bi-Hamiltonian system.
Conversely, if $Z$ is even-dimensional,
any Poisson bivector field $w$ (\ref{bi42})
is extended
to an appropriate symplectic structure $\Om$
as follows.

\begin{prop} \label{bi81} \mar{bi81}
The Poisson bivector field $w$ (\ref{bi42}) on a toroidal domain
$U'$ in Theorem \ref{bi92} is extended to a symplectic structure
$\Om$ on $U'$ such that integrals of motion $H_\la=J_\la$ remain in involution
and their Hamiltonian vector fields with respect to $w$ and $\Om$ coincide.
\end{prop}

\begin{proof}.
The Poisson bivector field $w$ (\ref{bi42}) on the foliated manifold 
$(U,\cF)$ defines a
leafwise symplectic form
$\Om_\cF$ (\ref{bi102}).
Restricted to the toroidal domain $U'$ in Theorem \ref{bi92} where
coordinates $J_\la$ have trivial transition functions, the exact 
sequence (\ref{bi121})
admits the splitting
\be
\zeta^*: T\cF^*\to T^*U', \qquad \zeta^*(\ol dJ_\m)= dJ_\m, \qquad 
\zeta^*(\ol dy^\m)= dy^\m
\ee
such that $\zeta^*\circ\Om_\cF$ is a presymplectic form on $U'$.
Let $\Om_Z=\Om_{AB}(z^C) dz^A\w dz^B$ be also a presymplectic form on 
$U'$. It always exist.
Then, $\Om=\zeta^*\circ\Om_\cF + \Om_Z$
  is a desired symplectic form on $U'$.
\end{proof}

\bigskip
\bigskip

\noindent
{\bf V. KAM THEOREM FOR PARTIALLY INTEGRABLE SYSTEMS}
\bigskip

Let $\{\cH_i\}$, $i=1,\ldots,k$, be a partially integrable system on a
$2n$-dimensional symplectic manifold $(Z,\Om)$.
Let $M$ be its regular connected
compact invariant manifold  which
admits an open neighbourhood satisfying Theorem \ref{dd3}.
In this case, Theorem \ref{dd3} comes to the
above mentioned Nekhoroshev theorem.
By virtue of this theorem,
there exists an
open neighbourhood of $M$ which is a trivial composite bundle
\mar{zz10}\beq
\pi:U=V\times W\times T^k\to V\times W\to V \label{zz10}
\eeq
(cf. (\ref{z10'})) over domains $W\subset \Bbb R^{2(n-k)}$ and 
$V\subset \Bbb R^k$. It
is provided with the partial
action-angle coordinates
$(I_i,z^A, \f^i)$, $i=1,\ldots,k$, $A=1,\ldots,2(n-k)$, such that
the symplectic form $\Om$ on $U$ reads
\mar{dd26}\beq
\Om= dI_i\w d\f^i +\Om_{AB}(I_j,z^C) dz^A\w dz^B +\Om_A^i(I_j,z^C) dI_i\w dz^A
\label{dd26}
\eeq
(cf. (\ref{d26})), while integrals of motion $H_i$ depend only on the action
coordinates $I_j$.

Note that, in accordance with part (III) of Theorem \ref{dd3},
one can always restrict $U$ to a Darboux coordinate
chart  provided with
coordinates $(I_i,p_s,q^s;\vf^i)$ such that the symplectic form $\Om$ 
(\ref{dd26})
takes the canonical form
\be
\Om= dI_i\w d\vf^i + dp_s\w dq^s.
\ee
Then, the PIS $\{H_i\}$ on this
chart can be extended to a completely
integrable system, e.g., $\{H_i,p_s\}$, but its invariant manifolds fail
to be compact. Therefore, this is not the case of
the KAM theorem.

Let $\cH(I_j)$ be a Hamiltonian of a PIS on $U$ (\ref{zz10}).
Its Hamiltonian vector field
\mar{k7}\beq
\xi=\dr^i\cH(I_\j)\dr_i \label{k7}
\eeq
with respect to the symplectic form $\Om$ (\ref{dd26})
yields the Hamilton equation
\mar{k65}\beq
\dot I_i=0, \qquad \dot z^A=0, \qquad \dot\f^i=\dr^i\cH(I_j) \label{k65}
\eeq
on $U$. Let us consider perturbations
\mar{k12'}\beq
\cH'=\cH+\cH_1(I_j,z^A,\f^j). \label{k12'}
\eeq
We assume the following.
(i) The Hamiltonian $\cH$ and its perturbations (\ref{k12'}) are real
analytic, although generalizations to the case of infinite and
finite order of differentiability are possible.$^{10,22}$
(ii) The Hamiltonian $\cH$ is nondegenerate, i.e., the
frequency map
\be
\om:V\times W\ni (I_j,z^A)\mapsto \xi^i(I_j)\in \Bbb R^k
\ee
is of rank $k$.

Note that $\om(V\times W)\subset \Bbb R^k$ is open
and bounded. As usual, given $\g>0$, let
\be
\Om_\g=\{ \om\in \Bbb R^k \,:\, |\om^ia_i|\geq \g(\op\sum^k_{j=1}
|a_j|)^{-k-1}, \quad \forall a\in\Bbb Z^k\setminus 0\}
\ee
denote the Cantor set of nonresonant frequences. The complement
of $\Om_\g\cap \om(V\times W)$ in $\om(V\times W)$ is dense and
open, but its relative Lebesgue measure tends to zero with $\g$.
Let us denote $\G_\g=\om^{-1}(\Om_\g)$, also called the Cantor set.

A problem is that the Hamiltonian vector field of the perturbed Hamiltonian
(\ref{k12'}) with respect to the symplectic form $\Om$
(\ref{dd26}) leads to the Hamilton equation $\dot z^A\neq 0$ and,
therefore, no torus (\ref{k65}) persists.

To overcome this difficulty, let us provide the toroidal domain
$U$ (\ref{zz10}) with the degenerate Poisson structure given by the
Poisson bivector field
\mar{k0}\beq
w=\dr^i\w \dr_i \label{k0}
\eeq
of rank $2k$. It is readily observed that, relative to $w$,
all integrals of motion of the original PIS $(\Om,\{H_i\})$ remain
in involution and, moreover, they possess the same Hamiltonian vector
fields $\vt_i$. In particular, a Hamiltonian $\cH$ with respect to the
Poisson structure (\ref{k0}) leads to the same Hamilton
equation (\ref{k65}). Thus, we can think of the pair $(w,\{H_i\})$
as being a PIS on the Poisson manifold $(U,w)$. The key point is that,
with respect to the Poisson bivector field $w$ (\ref{k0}),
the Hamiltonian vector
field of the perturbed Hamiltonian $\cH'$ (\ref{k12'}) is
\mar{bi51'}\beq
\xi'=\dr^i\cH'\dr_i -\dr_i\cH'\dr^i, \label{bi51'}
\eeq
and the corresponding first order dynamic equation on $U$ reads
\mar{k11}\beq
\dot I_i=-\dr_i\cH'(I_j,z^B,\f^j), \qquad \dot z^A=0,
\qquad \dot \f^i=\dr^i\cH'(I_j,z^B,\f^j). \label{k11}
\eeq
This is a Hamilton equation with respect to the Poisson structure
$w$ (\ref{k0}), but is not so relative to the original symplectic form
$\Om$. Since $\dot z^A=0$ and the toroidal domain $U$ (\ref{zz10})
is a trivial bundle over $W$, one can think of the dynamic equation
(\ref{k11}) as being
a perturbation of the dynamic equation (\ref{k65}) depending on parameters
$z^A$. Furthermore, the Poisson manifold $(U,w)$ is the
product of symplectic manifold $(V\times T^k,\Om')$
with the symplectic form
\mar{k4}\beq
\Om'=dI_i\w d\f^i \label{k4}
\eeq
and the Poisson manifold $(W,w=0)$ with the zero Poisson structure.
Therefore, the equation (\ref{k11}) can be seen as a Hamilton equation
on the symplectic manifold $(V\times T^k,\Om')$ depending on parameters.
Then, one can apply the conditions of quasi-periodic stability
of symplectic Hamiltonian systems depending on parameters$^{10}$
to the perturbation (\ref{k11}).

In a more general setting, these conditions can be formulated as follows.
Let $(w,\{H_i\})$, $i=2,\ldots,k$, be a PIS on a
regular Poisson manifold $(Z,w)$ of rank $2k$.
Let $M$ be its regular connected
compact invariant manifold, and let $U$ be its toroidal
neighbourhood $U$ (\ref{zz10}) in Theorem \ref{bi100}
provided with the partial
action-angle coordinates
$(I_i,z^A, \f^i)$  such that
the Poisson bivector $w$ on $U$ takes the canonical form (\ref{k0}).
The following result is a reformulation of that in Ref. [10] (Section 5c),
where $P=W$ is a
parameter space and $\si$ is the symplectic form (\ref{k4}) on $V\times T^k$.

\begin{theo} \label{cc14} \mar{cc14}
Given a torus
$\{0\}\times T^k$, let
\mar{k77}\beq
\xi=\xi^i(I_j,z^A)\dr_i\label{k77}
\eeq
(cf. (\ref{k7}))
be a real analytic Hamiltonian vector field whose frequency map
\be
\om: V\times W\ni (I_j,z^A)\mapsto \xi^i(I_j,z^A)\in \Bbb R^k
\ee
is of maximal rank at $\{0\}$. Then, there exists a neighbourhood 
$N_0\subset V\times W$ of $\{0\}$ such that,
for any real analytic Hamiltonian vector field
\be
\wt\xi=\wt\xi_i(I_j,z^A,\f^j)\dr^i + \wt\xi^i(I_j,z^A,\f^j)\dr_i
\ee
(cf. (\ref{bi51'})) sufficiently near $\xi$ (\ref{k77}) in the real 
analytic topology, the following holds.
Given the Cantor set $\G_\g\subset N_0$, there exists the $\wt\xi$-invariant
Cantor set $\wt\G\subset N_0\times T^k$ which is a 
$C^\infty$-near-identity diffeomorphic
image of $\G_\g\times T^k$.
\end{theo}

Theorem \ref{cc14} is an extension of the KAM theorem to PISs
on Poisson manifolds $(Z,w)$.
Given a PIS $(\Om,\{H_i\})$ on a symplectic manifold $(Z,\Om)$,
Theorem \ref{cc14} enables one to obtain its perturbations (\ref{bi51'})
possessing a large number of invariant tori, though these perturbations are
not Hamiltonian.

\bigskip
\bigskip

\noindent
{\bf VI. APPENDIX A}
\bigskip

{\it Proof of Theorem \ref{bi100}:}
First, let us employ Theorem \ref{bi92} and restrict $U$ to the toroidal
domain, say  $U$ again, equipped with coordinates $(J_\la,z^A,y^\la)$
such that the Poisson bivector field $w$ takes the form (\ref{bi42}) 
and the algebra
$\cA$ is generated by the Hamiltonian vector fields
$\vt_\la$ (\ref{bi93}) of $k$ independent functions
$H_\la=J_\la$ in involution. Let us choose these vector fields as new
generators of the group $G$ and return to Theorem \ref{bi0}. In accordance with
this theorem, there exists a toroidal domain $U'\subset U$ provided
with another trivialization $U'\to N'\subset N$ in toroidal cylinders
$\Bbb R^{k-m}\times T^m$ and endowed with bundle coordinates
$(J_\la,z^A,y'^\la)$ such that
the vector fields $\vt_\la$ (\ref{bi93}) take the form (\ref{ww25}).
For the sake of simplicity, let $U'$, $N'$ and $y'$ be denoted
$U$, $N$ and $y=(t^a,\vf^i)$ again. Herewith, the Poisson bivector field $w$
is given by the expression (\ref{bi42}) with new coefficients.

Let
$w^\sh: T^*U\to TU$ be the corresponding bundle homomorphism,
and let $T\cF^*\to U$ denote the dual of the characteristic distribution
$T\cF\to U$. We have the exact sequences
\mar{bi120,1}\ben
&& 0\to T\cF \ar^{i_\cF} TU \ar TU/T\cF\to 0, \label{bi120}\\
&& 0\to {\rm Ann}\,T\cF\ar T^*U\ar^{i^*_\cF} T\cF^* \to 0. \label{bi121}
\een
The bundle homomorphism $w^\sh$
factorizes in a unique fashion
\be
w^\sh:
T^*U\ar^{i^*_\cF} T\cF^*\ar^{w^\sh_\cF}
T\cF\ar^{i_\cF} TU
\ee
through the bundle isomorphism
\be
w_\cF^\sh: T\cF^*\to T\cF,  \qquad
w^\sh_\cF:\al\mapsto -w(x)\lfloor \al.
\ee
Then, the inverse isomorphisms
$w_\cF^\fl : T\cF\to T\cF^*$ provides the foliated manifold
$(U,\cF)$ with the
leafwise symplectic form
\mar{bi102,'}\ben
&& \Om_\cF=\Om^{\m\n}(J_\la,z^A,t^a) \ol dJ_\m\w \ol dJ_\n +
\Om_\m^\n(J_\la,z^A) \ol dJ_\n\w \ol dy^\m, \label{bi102}\\
&& \Om_\m^\al w^\m_\bt=\dl^\al_\bt, \qquad 
\Om^{\al\bt}=-\Om^\al_\m\Om^\bt_\n w^{\m\n},
\label{bi102'}
\een
where $\{\ol dJ_\m, \ol dy^\m\}$ is the dual of the
basis $\{\dr^\m,\dr_\m\}$ for the characteristic distribution $T\cF$.
Recall that leafwise (or tangential) exterior forms are defined as sections of
the exterior bundle $\w T\cF^*\to U$, while the
leafwise exterior differential $\ol d$
acts on them by the law
\be
\ol d\psi=\ol dJ_\la\w\dr^\la\psi +\ol dy^\la\w\dr_\la\psi
\ee
(see, e.g., Ref. [23,24]). The leafwise symplectic form $\Om_\cF$
is nondegenerate and $\ol d$-closed, i.e., $\ol d\Om_\cF=0$.
Let us show that it is $\ol d$-exact.

Let $F$ be a leaf of the foliation $\cF$ of $U$. There is a homomorphism
of the de Rham cohomology $H^*(U)$ of $U$ to the de Rham
cohomology of $H^*(F)$ of $F$. One can show that this homomorphism factorizes
through the leafwise cohomology$^{24}$
\mar{lmp11}\beq
H^*(U)\to H^*_\cF(U)\to H^*(F). \label{lmp11}
\eeq
Since $N$ is a domain of an adapted coordinate chart of
the foliation $\cF_N$, the foliation $\cF_N$ of $N$ is a trivial
fiber bundle $N=V\times W\to W$. Since $\cF$ is the pull-back onto $U$
of the foliation $\cF_N$ of $N$, it is also a trivial fiber bundle
\mar{bi103}\beq
U=V\times W\times (\Bbb R^{k-m}\times T^m) \to W \label{bi103}
\eeq
over a domain $W\subset \Bbb R^{\di Z-2k}$. It follows that
\be
H^*(U)=H^*(T^m)=H^*_\cF(U).
\ee
Then, the closed leafwise two-form $\Om_\cF$ (\ref{bi102})
is exact due to the absence of the term $\Om_{\m\n}dy^\m\w dy^\nu$.
Moreover, $\Om_\cF=\ol d\Xi$ where $\Xi$ reads
\be
\Xi=\Xi^\al(J_\la,z^A,y^\la)\ol dJ_\al + \Xi_i(J_\la,z^A)\ol d\vf^i
\ee
up to a $\ol d$-exact leafwise form.

The Hamiltonian vector fields $\vt_\la=\vt_\la^\m\dr_\m$ (\ref{ww25})
obey the relation
\mar{ww22'}\beq
\vt_\la\rfloor\Om_\cF=-\ol dJ_\la, \qquad \Om^\al_\bt 
\vt^\bt_\la=\dl^\al_\la, \label{ww22'}
\eeq
which falls into the following conditions
\mar{bi110,1}\ben
&& \Om^\la_i=\dr^\la\Xi_i-\dr_i\Xi^\la, \label{bi110} \\
&& \Om^\la_a=-\dr_a\Xi^\la=\dl^\la_a. \label{bi111}
\een
The first of the relations (\ref{bi102'})  shows that $\Om^\al_\bt$ 
is a nondegenerate matrix
independent of coordinates $y^\la$. Then, the condition (\ref{bi110}) implies
that $\dr_i\Xi^\la$ are
independent of $\vf^i$, and so are $\Xi^\la$ since $\vf^i$ are cyclic 
coordinates. Hence,
\mar{bi112,3}\ben
&&\Om^\la_i=\dr^\la\Xi_i, \label{bi112}\\
&& \dr_i\rfloor\Om_\cF=-\ol d\Xi_i. \label{bi113}
\een
Let us introduce new coordinates $I_a=J_a$, $I_i=\Xi_i(J_\la)$. By 
virtue of the
equalities (\ref{bi111}) and (\ref{bi112}), the Jacobian of this coordinate
transformation is regular. The relation
(\ref{bi113}) shows that $\dr_i$ are Hamiltonian vector fields of the functions
$H_i=I_i$. Consequently, we can choose vector fields
$\dr_\la$ as generators of the algebra $\cA$.
One obtains from the equality (\ref{bi111}) that
$\Xi^a=-t^a+E^a(J_\la,z^A)$
and $\Xi^i$ are independent of $t^a$. Then, the leafwise Liouville 
form $\Xi$ reads
\be
\Xi=(-t^a+E^a(I_\la,z^A))\ol dI_a + E^i(I_\la, z^A)\ol dI_i + I_i \ol d\vf^i.
\ee
The coordinate shifts
\be
x^a=-t^a+E^a(I_\la,z^A), \qquad \f^i=\vf^i-E^j(I_\la,z^A)
\ee
bring the leafwise form $\Om_\cF$ (\ref{bi102}) into the canonical form
\be
\Om_\cF= \ol dI_a\w \ol d x^a + \ol dI_i\w \ol d\f^i
\ee
which ensures the canonical form (\ref{bi101}) of the Poisson 
bivector field $w$.

\bigskip
\bigskip

\noindent
{\bf VII. APPENDIX B}
\bigskip

{\it Proof of Theorem \ref{dd3}:} (I). See
the proof of part (I) of Theorem \ref{bi0}.

(II). One can specify the coordinates on the
base $N$ of the trivial bundle $U\to N$ as follows. Let us
consider the morphism
\mar{d11}\beq
\pi'=\op\times^\la H_\la: U\to V \label{d11}
\eeq
of $U$ onto a domain $V\subset \Bbb R^k$. It is of constant rank and,
consequently, is a fibered manifold. The fibration $\pi'$ factorizes as
\be
\pi': U\ar^\pi N\ar^{\pi''} V
\ee
through the fiber bundle $\pi$. The map $\pi''=\pi'\circ\si$
is also a fibered manifold. One can always restrict
the domain $N$ to a chart of the fibered manifold $\pi''$. Then,
$N\to\pi''(N)=V$ is a trivial bundle, and so is $U\to V$.
Thus, we have the composite fibration
\mar{z10'}\beq
U=V\times W\times (\Bbb R^{k-m}\times T^m)\to V\times W\to V \label{z10'}
\eeq
Let us provide its base $V$ with the coordinates $(J_\la)$ such that
$J_\la(u)=H_\la(u)$, $u\in U$. Then $N$ can be equipped with the
bundle coordinates $(J_\la, z^A)$, $A=1,\ldots, 2(n-k)$, and $(J_\la, z^A,
t^a,\vf^i)$ are coordinates on $U$ (\ref{z10'}).
Since fibers of $U\to N$ are isotropic,
the symplectic form $\Om$ on $U$
relative to the coordinates $(J_\la,z^A,y^\la)$ reads
\mar{d23}\beq
\Om=\Om^{\al\bt}dJ_\al\w dJ_\bt + \Om^\al_\bt dJ_\al\w dy^\bt +
\Om_{AB}dz^A\w dz^B +\Om_A^\la dJ_\la\w dz^A +
  \Om_{A\bt} dz^A\w dy^\bt. \label{d23}
\eeq
The Hamiltonian vector fields $\vt_\la=\vt_\la^\m\dr_\m$ (\ref{ww25})
obey the relations $\vt_\la\rfloor\Om=-dJ_\la$, which give the 
coordinate conditions
\mar{ww22}\beq
  \Om^\al_\bt \vt^\bt_\la=\dl^\al_\la, \qquad \Om_{A\bt}\vt^\bt_\la=0. 
\label{ww22}
\eeq
The first of them shows that $\Om^\al_\bt$ is a nondegenerate matrix
independent of coordinates $y^\la$.
Then, the second one implies $\Om_{A\bt}=0$.

By virtue of the well-known K\"unneth formula for the de Rham cohomology
of manifold products, the closed form $\Om$ (\ref{d23}) is exact, i.e.,
$\Om=d\Xi$ where the Liouville form $\Xi$ is
\be
\Xi=\Xi^\al(J_\la,z^B,y^\la)dJ_\al + \Xi_i(J_\la,z^B) d\vf^i
+\Xi_A(J_\la,z^B,y^\la)dz^A.
\ee
Since $\Xi_a=0$ and $\Xi_i$ are independent of $\vf^i$, it follows from
the relations
\be
\Om_{A\bt}=\dr_A\Xi_\bt-\dr_\bt\Xi_A=0
\ee
that $\Xi_A$ are independent of coordinates $t^a$ and are
at most affine in $\vf^i$. Since $\vf^i$ are cyclic coordinates,
  $\Xi_A$ are independent of $\vf^i$. Hence,
$\Xi_i$ are  independent of coordinates $z^A$,
and the Liouville form reads
\mar{ac2}\beq
\Xi=\Xi^\al(J_\la,z^B,y^\la)dJ_\al + \Xi_i(J_\la) d\vf^i
+\Xi_A(J_\la,z^B)dz^A. \label{ac2}
\eeq
Because entries $\Om^\al_\bt$ of $d\Xi=\Om$ are independent of 
$y^\la$, we obtain the following.

(i) $\Om^\la_i=\dr^\la\Xi_i-\dr_i\Xi^\la$. Consequently, $\dr_i\Xi^\la$ are
independent of $\vf^i$, and so are $\Xi^\la$ since $\vf^i$ are cyclic 
coordinates. Hence,
$\Om^\la_i=\dr^\la\Xi_i$ and $\dr_i\rfloor\Om=-d\Xi_i$. A glance at 
the last equality
shows that $\dr_i$ are Hamiltonian vector fields. It follows that, 
from the beginning,
one can separate $m$ integrals of motion, say $H_i$ again, whose Hamiltonian
vector fields are tangent to invariant tori. In this case, the matrix 
$B$ in the expressions
(\ref{d6}) and (\ref{ww25}) vanishes, and the Hamiltonian vector fields
$\vt_\la$ (\ref{ww25}) read
\mar{ww25'}\beq
\vt_a=\dr_a, \qquad \vt_i=(C^{-1})_i^k\dr_k. \label{ww25'}
\eeq
Moreover, the coordinates $t^a$ are exactly the flow parameters $s^a$.
Substituting the expressions (\ref{ww25'}) into the first condition 
(\ref{ww22}), we obtain
\be
\Om=\Om^{\al\bt}dJ_\al\w dJ_\bt +dJ_a\w ds^a + C^i_k dJ_i\w d\vf^k +
\Om_{AB}dz^A\w dz^B +\Om_A^\la dJ_\la\w dz^A.
\ee
It follows that $\Xi_i$ are independent of $J_a$, and so are 
$C^k_i=\dr^k\Xi_i$.

(ii) $\Om^\la_a=-\dr_a\Xi^\la=\dl^\la_a$. Hence, $\Xi^a=-s^a+E^a(J_\la)$
and $\Xi^i$ are independent of $s^a$.

In view of items (i) -- (ii), the Liouville form $\Xi$ (\ref{ac2}) reads
\be
\Xi=(-s^a+E^a(J_\la,z^B))dJ_a + E^i(J_\la,z^B) dJ_i + \Xi_i(J_j) d\vf^i +
\Xi_A(J_\la,z^B)dz^A.
\ee
Since the matrix $\dr^k\Xi_i$ is nondegenerate,
we can perform the coordinate transformation $I_a=J_a$, 
$I_i=\Xi_i(J_j)$ together
with the coordinate shifts
\be
x^a=-s^a+E^a(J_\la,z^B), \qquad \f^i=\vf^i-E^j(J_\la,z^B)\frac{\dr 
J_j}{\dr I_i}.
\ee
These transformations bring $\Om$ into the form (\ref{d26}).

(III). Since functions $I_\la$ are in involution and their Hamiltonian
vector fields $\dr_\la$ mutually commute, a
point $z\in M$ has an open neigbourhood $Q\times O_z$, $O_z\in
\Bbb R^{k-m}\times T^m,$  endowed with the
Darboux coordinates $(I_\la,p_s,q^s,\wt y^\la)$
such that the symplectic form $\Om$ (\ref{d26}) is given by the expression
(\ref{dd12}).
Here, $\wt y^\la(I_\la,z^A,y^\al)$ are local functions whose Hamiltonian
vector fields are $\dr^\la$. They take
the form
\mar{dd11}\beq
\wt y^\la=y^\la + f^\la(I_\la,z^A). \label{dd11}
\eeq
With the group $G$, one can extend these functions to
the open neighbourhood
\be
\wt U=Q\times \Bbb R^{k-m}\times T^m
\ee
of $M$ by the law
\be
\wt y^\la(I_\la,z^A,G(y)^\al)= G(y)^\la + f^\la(I_\la,z^A).
\ee
Substituting the functions (\ref{dd11}) on $\wt U$
into the expression (\ref{d26}),
one brings the symplectic form $\Om$ into the canonical form
(\ref{dd12}) on $\wt U$.

\end{document}